\documentclass[11pt, a4paper]{article}

\usepackage{amsmath, amssymb, amsthm}
\usepackage{graphicx}
\usepackage{geometry}
\usepackage{listings}
\usepackage{xcolor}
\usepackage{hyperref}
\usepackage{booktabs}
\usepackage{algorithm}
\usepackage{algpseudocode}
\usepackage{authblk}
\usepackage{float}
\usepackage{subcaption}

\usepackage{pgfplots}
\pgfplotsset{compat=1.18} 

\usepackage{tikz}
\usetikzlibrary{positioning, fit, arrows.meta, calc}

\title{Mesh-Free Numerical Approximation of the Biharmonic Equation via Optimized Kolmogorov–Arnold Neural Networks}

\author{B. Veena S. N. Rao \thanks{bv.rao@tamucc.edu}}
\affil{Department of Mathematics \& Statistics,\\
Texas A\&M University-Corpus Christi, \\
6300 Ocean Dr., 
Corpus Christi, Texas 78412-5825, USA}
\date{}

\begin{document}
\maketitle

\begin{abstract}
A mesh-free numerical framework based on \textit{Kolmogorov--Arnold Physics-Informed Neural Networks} (KAN-PINNs) is developed for the approximation of fourth-order elliptic boundary value problems, with specific application to the biharmonic equation governing thin plate deflection. Unlike conventional Multi-Layer Perceptrons relying on fixed nodal activations, learnable univariate functions parameterized via radial basis functions are deployed on network edges, while high-order differentiability is preserved through hyperbolic tangent activation mappings. The severe numerical stiffness inherent to fourth-order differential operators and fully clamped boundary conditions is addressed through a direct normalized residual formulation coupled with a hybrid, multi-stage AdamW-to-L-BFGS optimization pipeline. An automated \textit{24/7 hill-climbing search protocol} is implemented to systematically calibrate boundary penalty weights and optimization schedules. When evaluated on a smooth manufactured benchmark on the unit square, an error reduction factor exceeding $150\times$ is attained over successive iterations, culminating in a final relative $L_2$ error of $1.593 \times 10^{-5}$ ($0.00159\%$) and a training loss of $3.065 \times 10^{-6}$ utilizing only $7,801$ trainable parameters. These results demonstrate that high-order PDE problems can be accurately resolved using compact, mesh-free KAN-PINNs architecture without requiring auxiliary variable transformations or discrete mesh generation.

\vspace{0.5em}
\end{abstract}

\section{Introduction}
\label{sec:introduction}

Partial differential equations (PDEs) provide a fundamental mathematical framework for describing a wide range of physical, engineering, and scientific phenomena, including heat and mass transfer, wave propagation, fluid flow, elasticity, diffusion, reaction processes, and many other multiscale systems. Although analytical solutions are available for a limited class of idealized problems, practical PDE models frequently involve complex geometries, heterogeneous material properties, nonlinearities, high-dimensional parameter spaces, and mixed boundary conditions. Consequently, reliable numerical approximation of PDE solutions remains an enduring cornerstone in computational science and engineering.

Among the most established approaches for the numerical solution of PDEs are finite element and Galerkin methods. The finite element method (FEM) provides a systematic framework for converting continuous boundary-value problems into finite-dimensional algebraic systems through variational formulations and localized basis functions. The mathematical foundations of FEM, including Sobolev spaces, variational formulations, approximation theory, stability, and convergence, have been extensively developed in the classical literature~\cite{brenner2008mathematical, ciarlet2002finite}. In particular, the formulation of FEM as a Galerkin approximation provides a natural connection between the underlying continuous problem and its discrete numerical representation. Strang and Fix~\cite{strang1973analysis} established many of the fundamental concepts underlying finite element approximation, including the relationship between variational formulations, Galerkin methods, approximation spaces, and convergence.

Spectral methods represent another important class of classical numerical techniques for PDEs. Unlike the locally supported polynomial basis functions commonly used in FEM, spectral methods employ globally supported, high-order basis functions to approximate the solution. For sufficiently smooth problems, this global approximation can lead to very high convergence rates and accurate numerical solutions. The theoretical and computational foundations of spectral methods, including polynomial approximation, Galerkin formulations, stability, convergence, and numerical solution procedures, are discussed extensively by Canuto et al.~\cite{canuto2006spectral}. These methods provide important high-accuracy reference baselines against which newer numerical approaches can be assessed.

For fourth-order boundary value problems, such as those governed by the biharmonic operator arising in plate bending and structural mechanics, classical discretization presents additional mathematical and computational complexities. Finite difference schemes~\cite{Ehrlich1971, BenArtzi2009}, mixed finite element formulations~\cite{Monk1987, CiarletRaviart1974}, and compact finite difference methods~\cite{Pan2025} require careful treatment of high-order derivative coupling and multiple boundary conditions. While these classical techniques are mathematically mature, their implementation on complex domains or higher dimensions can involve substantial mesh generation overhead and large sparse algebraic systems.

These computational demands, alongside the increasing availability of physical and experimental data, have motivated the development of numerical approaches capable of combining governing differential equations with data-driven learning. These developments have contributed to the rapid emergence of scientific machine learning and, in particular, neural-network-based methods for solving PDEs. An influential development in this direction is the physics-informed neural network (PINN). Raissi et al.~\cite{raissi2019physicsinformed} introduced a general framework in which neural networks are trained while simultaneously enforcing the governing differential equations and associated initial or boundary conditions. Instead of constructing a conventional mesh and explicitly assembling a discrete PDE operator, a neural network is used as a global approximation of the unknown solution, while automatic differentiation is employed to evaluate the derivatives appearing in the governing equations. This formulation enables the integration of physical constraints into the learning process and has been applied to both forward and inverse PDE problems. Beyond multi-dimensional PDE systems, feedforward neural network architectures have been successfully adapted to approximate solutions for a wide spectrum of linear and nonlinear ordinary, singular, and delay differential equations~\cite{venkatachalapathy2023feedforward, venkatachalapathy2023deep, mallikarjunaiah2023deep}. The broader role of physics-informed machine learning in scientific computing has subsequently been reviewed by Karniadakis et al.~\cite{karniadakis2021physics}.

The development of neural PDE solvers has also revealed important limitations associated with the optimization of physics-informed loss functions. In particular, Wang et al.~\cite{wang2021understanding} investigated gradient-flow pathologies in PINNs and demonstrated that the different components of the loss function can produce highly unbalanced optimization dynamics. Such behavior can make training difficult, particularly when differential operators introduce disparate spatial scales or higher-order differentiability requirements. Krishnapriyan et al.~\cite{krishnapriyan2021characterizing} further demonstrated that PINN formulations can exhibit failure modes for challenging PDEs, showing that the difficulties are not necessarily caused by insufficient neural-network expressivity but can instead arise from the structure and optimization of the physics-informed objective. For higher-order problems like the biharmonic equation, these optimization hurdles are amplified due to the repeated differentiation of network parameters. Consequently, recent studies have explored specialized PINN architectures and training strategies specifically tailored for fourth-order systems~\cite{DwivediSrinivasan2020, Vahab2022, liu2023gradient, MishraKhan2026}.

Alternative neural formulations have also been established to address these optimization dynamics. The Deep Ritz method proposed by E and Yu~\cite{e2018deep} represents the solution using a neural network and minimizes the corresponding variational or energy functional rather than pointwise differential residuals. This approach establishes a direct connection between variational numerical methods and neural-network approximation. Similarly, the framework of Galerkin Neural Networks developed by Ainsworth and Dong~\cite{ainsworth2021galerkin} constructs finite-dimensional approximation spaces whose basis functions are generated by neural networks and subsequently uses these spaces within a Galerkin formulation, preserving classical error-control concepts while maintaining neural flexibility.

More recently, the introduction of Kolmogorov-Arnold Networks (KANs) by Liu et al.~\cite{liu2024kan} has provided an alternative neural-network architecture for scientific computing. In contrast to standard multilayer perceptrons (MLPs), in which fixed nonlinear activation functions are associated with neurons and trainable scalar weights connect layers, KANs employ learnable univariate functions parameterized as B-splines on the edges of the network. This architectural difference provides an alternative mechanism for representing complex functional relationships and has motivated investigations of KANs in scientific and PDE-related applications.

The combination of KAN architectures with physics-informed learning has subsequently led to \textit{Kolmogorov-Arnold-Informed Neural Networks} (KINNs). Wang et al.~\cite{wang2025kinn} introduced KINN as a physics-informed framework in which KANs are incorporated into strong-form, energy-form, and inverse-form PDE formulations. This provides a direct connection between the architectural flexibility of KANs and the physics-constrained learning paradigm of PINNs. Furthermore, specialized KAN paradigms—such as physically consistent and interpretable SLE-KAN frameworks—have proven highly effective in advanced constitutive modeling of soft materials~\cite{pati2026sle}, highlighting the versatility of edge-based neural parameterizations.

These developments establish a natural methodological spectrum extending from classical discretization techniques to modern neural PDE solvers. At one end of this spectrum, FEM and conventional Galerkin methods rely on carefully constructed finite-dimensional approximation spaces with well-established mathematical theories of stability and convergence~\cite{brenner2008mathematical, strang1973analysis, ciarlet2002finite}. Spectral methods employ global high-order approximations and achieve rapid convergence for smooth solutions~\cite{canuto2006spectral}. To bridge these paradigms, detailed comparative assessments between physics-informed deep learning and classical finite element methods have been conducted across complex transport systems, including Darcy--Brinkman--Forchheimer porous media flow and electrohydrodynamic formulations~\cite{martinez2024approximation, martinez2026numerical}, illuminating the relative strengths and computational trade-offs of mesh-free neural approximations. Neural approaches replace or augment these conventional approximation spaces with trainable function representations. PINNs enforce differential equations through residual-based loss functions~\cite{raissi2019physicsinformed}, while Deep Ritz and Galerkin neural approaches establish direct connections with variational numerical analysis~\cite{e2018deep, ainsworth2021galerkin}. Finally, KAN and KINN introduce a distinct neural architecture and integrate it with physics-informed frameworks~\cite{liu2024kan, wang2025kinn}.

Building upon this broader methodological progression, the objective of this study is to formulate, implement, and rigorously evaluate a mesh-free KAN-PINNs framework specifically designed for approximating solutions to the biharmonic equation under fully clamped boundary conditions. By integrating a direct fourth-order residual formulation, smooth activation mappings, and an automated multi-stage optimization pipeline, the proposed approach successfully overcomes the numerical stiffness characteristic of high-order plate bending problems. The remainder of this paper is organized as follows: Section~\ref{sec:math_foundations} establishes the mathematical formulation, weak analysis, and governing residual equations. Section~\ref{sec:architecture} details the KAN layer mechanics, architectural topology, and approximation guarantees. Section~\ref{sec:optimization} outlines the hybrid optimization strategy and the automated hill-climbing search protocol. Comprehensive numerical results, pointwise error distributions, and convergence trajectories are presented and analyzed in Section~\ref{sec:results}, followed by concluding remarks in Section~\ref{sec:conclusion}.


\section{Mathematical foundations}
\label{sec:math_foundations}

This section establishes the mathematical formulation of the mesh-free Kolmogorov--Arnold Physics-Informed Neural Network (KAN-PINN) for the solution of the biharmonic equation associated with the transverse deflection of a thin, clamped elastic plate. We first formulate the governing boundary-value problem, discuss its weak formulation and well-posedness, and subsequently introduce the KAN-PINN approximation and the corresponding physics-informed loss functional.

\subsection{The biharmonic equation}
\label{subsec:governing_equation}

Let $\Omega = (0,1)^2 \subset \mathbb{R}^2$ denote the bounded domain occupied by a thin elastic plate, and let $\partial\Omega$ denote its boundary. Under the assumptions of small transverse deflection, linear elasticity, and Kirchhoff--Love plate theory, the transverse displacement $u:\overline{\Omega} \rightarrow \mathbb{R}$ satisfies the biharmonic boundary-value problem
\begin{equation}
    \Delta^2 u = f
    \qquad \text{in } \Omega,
    \label{eq:biharmonic}
\end{equation}
where $f:\Omega\rightarrow\mathbb{R}$ denotes the prescribed transverse load and $\Delta^2 u = \Delta(\Delta u)$ is the biharmonic operator. In two spatial dimensions, the operator can be written explicitly as
\begin{equation}
    \Delta^2 u
    =
    \frac{\partial^4 u}{\partial x^4}
    +
    2\frac{\partial^4 u}{\partial x^2\partial y^2}
    +
    \frac{\partial^4 u}{\partial y^4}.
    \label{eq:biharmonic_expanded}
\end{equation}

For a fully clamped plate, both the transverse displacement and the normal slope are prescribed along the boundary. Thus, the boundary conditions are
\begin{equation}
    u = g
    \qquad \text{on } \partial\Omega,
    \label{eq:dirichlet_bc}
\end{equation}
and
\begin{equation}
    \frac{\partial u}{\partial n}
    =
    \nabla u\cdot n
    =
    j
    \qquad \text{on } \partial\Omega,
    \label{eq:neumann_bc}
\end{equation}
where $n$ denotes the outward unit normal vector to $\partial\Omega$, while $g$ and $j$ represent the prescribed displacement and normal slope, respectively.

Consequently, the complete clamped biharmonic problem can be stated as
\begin{equation}
\begin{aligned}
    \Delta^2 u &= f && \text{in } \Omega,\\
    u &= g && \text{on } \partial\Omega,\\
    \frac{\partial u}{\partial n} &= j
    && \text{on } \partial\Omega.
\end{aligned}
\label{eq:clamped_bvp}
\end{equation}
For a homogeneous clamped plate, the boundary data reduce to $g=0, j=0$, and hence
\begin{equation}
    u=0,
    \qquad
    \frac{\partial u}{\partial n}=0
    \qquad \text{on } \partial\Omega.
    \label{eq:homogeneous_clamped}
\end{equation}

\subsection{Weak formulation}
\label{subsec:weak_formulation}

Because the biharmonic operator is of fourth order, the natural energy space for the clamped plate problem is a subspace of $H^2(\Omega)$. For the homogeneous boundary-value problem, define
\begin{equation}
    H_0^2(\Omega)
    =
    \left\{
    v\in H^2(\Omega):
    v=0,\;
    \frac{\partial v}{\partial n}=0
    \text{ on }\partial\Omega
    \right\},
    \label{eq:H02}
\end{equation}
where the boundary conditions are understood in the sense of Sobolev traces.

Let $v\in H_0^2(\Omega)$ be an admissible test function. Multiplying \eqref{eq:biharmonic} by $v$ and integrating over $\Omega$ gives
\begin{equation}
    \int_{\Omega} (\Delta^2 u)v\,d\mathbf{x}
    =
    \int_{\Omega} fv\,d\mathbf{x}.
    \label{eq:weak_start}
\end{equation}

After integration by parts and using the clamped boundary conditions, the corresponding variational formulation can be expressed in terms of the Hessian as
\begin{equation}
    a(u,v)
    =
    \ell(v),
    \qquad
    \forall v\in H_0^2(\Omega),
    \label{eq:weak_form}
\end{equation}
where
\begin{equation}
    a(u,v)
    =
    \int_{\Omega}
    D^2u : D^2v\,d\mathbf{x},
    \label{eq:bilinear_form}
\end{equation}
and
\begin{equation}
    \ell(v)
    =
    \int_{\Omega}fv\,d\mathbf{x}.
    \label{eq:linear_functional}
\end{equation}

Here, $D^2u$ denotes the Hessian matrix of $u$, and the Frobenius inner product is defined by
\begin{equation}
    D^2u:D^2v
    =
    \sum_{i,j=1}^{2}
    \frac{\partial^2u}{\partial x_i\partial x_j}
    \frac{\partial^2v}{\partial x_i\partial x_j}.
\end{equation}

Equivalently, under the clamped boundary conditions, one may write the bilinear form in terms of the Laplacian as
\begin{equation}
    a(u,v)
    =
    \int_{\Omega}
    \Delta u\,\Delta v\,d\mathbf{x}.
    \label{eq:laplacian_bilinear}
\end{equation}

The weak formulation therefore seeks $u\in H_0^2(\Omega)$ such that
\begin{equation}
    \int_{\Omega}
    D^2u:D^2v\,d\mathbf{x}
    =
    \int_{\Omega}fv\,d\mathbf{x},
    \qquad
    \forall v\in H_0^2(\Omega).
    \label{eq:weak_biharmonic}
\end{equation}

\subsection{Existence and uniqueness of the solution}
\label{subsec:existence_uniqueness}

The variational formulation provides the appropriate framework for establishing the existence and uniqueness of the clamped plate solution. Assume, for example, that $f\in H^{-2}(\Omega)$, where $H^{-2}(\Omega)$ denotes the dual space of $H_0^2(\Omega)$. For the commonly considered case of a sufficiently regular load, one may in particular take $f\in L^2(\Omega)$.

The bilinear form $a(\cdot,\cdot)$ defined in \eqref{eq:bilinear_form} is continuous on $H_0^2(\Omega)$. Moreover, for the clamped boundary conditions, an appropriate second-order Poincar\'e inequality implies the coercivity relation
\begin{equation}
    a(v,v)
    =
    \int_{\Omega}|D^2v|^2\,d\mathbf{x}
    \geq
    C\|v\|_{H^2(\Omega)}^2,
    \qquad
    \forall v\in H_0^2(\Omega),
    \label{eq:coercivity}
\end{equation}
for some constant $C>0$ depending on $\Omega$.

The linear functional $\ell$ is continuous on $H_0^2(\Omega)$. Indeed, for $f\in L^2(\Omega)$,
\begin{equation}
    |\ell(v)|
    =
    \left|
    \int_{\Omega}fv\,d\mathbf{x}
    \right|
    \leq
    \|f\|_{L^2(\Omega)}
    \|v\|_{L^2(\Omega)}
    \leq
    C\|f\|_{L^2(\Omega)}
    \|v\|_{H^2(\Omega)}.
    \label{eq:linear_continuity}
\end{equation}

Consequently, the Lax--Milgram theorem guarantees that there exists a unique weak solution $u\in H_0^2(\Omega)$ satisfying \eqref{eq:weak_biharmonic}. In particular, the solution depends continuously on the forcing term, with the stability estimate
\begin{equation}
    \|u\|_{H^2(\Omega)}
    \leq
    C\|f\|_{H^{-2}(\Omega)}.
    \label{eq:stability}
\end{equation}

For sufficiently smooth forcing and boundary data, elliptic regularity yields additional regularity of the solution. In particular, under appropriate domain and compatibility assumptions, the weak solution may possess sufficient regularity to satisfy the biharmonic equation in the classical sense. For example, when the solution belongs to $H^4(\Omega)$, the fourth-order derivatives in \eqref{eq:biharmonic_expanded} are well-defined in the $L^2$ sense, and the strong formulation \eqref{eq:clamped_bvp} is recovered.

For the square domain considered here, regularity at the corners should be interpreted with care. The regularity of a fourth-order elliptic problem may depend on the compatibility of the boundary data and on corner singularities. Therefore, when classical $C^4(\overline{\Omega})$ regularity is assumed for the PINN formulation, it is understood that the forcing and boundary data satisfy the corresponding regularity and compatibility requirements.

\subsection{Nonhomogeneous clamped boundary conditions}
\label{subsec:nonhomogeneous}

For nonhomogeneous boundary data $(g,j)$, the solution can be decomposed as
\begin{equation}
    u=w+\widetilde{g},
    \label{eq:lifting_decomposition}
\end{equation}
where $\widetilde{g}\in H^2(\Omega)$ is a lifting function satisfying
\begin{equation}
    \widetilde{g}=g,
    \qquad
    \frac{\partial\widetilde{g}}{\partial n}=j
    \qquad
    \text{on }\partial\Omega,
    \label{eq:lifting_bc}
\end{equation}
and $w\in H_0^2(\Omega)$.

Substitution into the governing equation gives
\begin{equation}
    \Delta^2 w
    =
    f-\Delta^2\widetilde{g}
    \qquad\text{in }\Omega,
    \label{eq:lifted_pde}
\end{equation}
with homogeneous clamped boundary conditions
\begin{equation}
    w=0,
    \qquad
    \frac{\partial w}{\partial n}=0
    \qquad\text{on }\partial\Omega.
\end{equation}

Thus, the same existence and uniqueness arguments apply to the nonhomogeneous problem after the boundary data have been incorporated through an appropriate lifting function.

\subsection{Mesh-free KAN approximation}
\label{subsec:kan_approximation}

To approximate the unique solution of \eqref{eq:clamped_bvp}, we employ a Kolmogorov--Arnold Network (KAN) as a mesh-free function approximator. Let
\begin{equation}
    u_{\boldsymbol{\theta}}(x,y)
    =
    \mathcal{N}_{\boldsymbol{\theta}}(x,y)
\end{equation}
denote the KAN approximation to the exact displacement field, where $\boldsymbol{\theta}$ represents the trainable parameters of the network.

Unlike mesh-based discretization methods, the KAN-PINN does not require a predefined mesh, element connectivity, or nodal basis. Instead, the governing differential equation and boundary conditions are enforced at a collection of collocation points sampled directly from the computational domain.

Let
\begin{equation}
    \mathcal{X}_{\Omega}
    =
    \left\{
    \mathbf{x}_k=(x_k,y_k)
    \right\}_{k=1}^{N_{\Omega}}
    \subset\Omega
\end{equation}
denote the set of interior collocation points, and let
\begin{equation}
    \mathcal{X}_{\partial\Omega}
    =
    \left\{
    \mathbf{x}_k^{\,b}
    \right\}_{k=1}^{N_{\partial\Omega}}
    \subset\partial\Omega
\end{equation}
denote the set of boundary collocation points.

\subsection{Physics-Informed Residual}
\label{subsec:pde_residual}

The differential equation is imposed through the strong-form residual
\begin{equation}
    r_{\mathrm{PDE}}(\mathbf{x};\boldsymbol{\theta})
    =
    \Delta^2u_{\boldsymbol{\theta}}(\mathbf{x})
    -
    f(\mathbf{x}).
    \label{eq:pde_residual}
\end{equation}

For $\mathbf{x}=(x,y)$, this becomes
\begin{equation}
    r_{\mathrm{PDE}}(x,y;\boldsymbol{\theta})
    =
    \frac{\partial^4u_{\boldsymbol{\theta}}}{\partial x^4}
    +
    2\frac{\partial^4u_{\boldsymbol{\theta}}}
    {\partial x^2\partial y^2}
    +
    \frac{\partial^4u_{\boldsymbol{\theta}}}{\partial y^4}
    -
    f(x,y).
    \label{eq:expanded_residual}
\end{equation}

The derivatives appearing in this expression are evaluated automatically using the differentiable structure of the KAN and automatic differentiation.

\subsection{Boundary residuals}
\label{subsec:boundary_residuals}

The clamped boundary conditions are imposed by defining two boundary residuals. The displacement residual is
\begin{equation}
    r_u(\mathbf{x};\boldsymbol{\theta})
    =
    u_{\boldsymbol{\theta}}(\mathbf{x})
    -
    g(\mathbf{x}),
    \qquad
    \mathbf{x}\in\partial\Omega,
    \label{eq:displacement_residual}
\end{equation}
while the normal-slope residual is
\begin{equation}
    r_n(\mathbf{x};\boldsymbol{\theta})
    =
    \nabla u_{\boldsymbol{\theta}}(\mathbf{x})
    \cdot n(\mathbf{x})
    -
    j(\mathbf{x}),
    \qquad
    \mathbf{x}\in\partial\Omega.
    \label{eq:normal_residual}
\end{equation}

\subsection{KAN-PINN loss functional}
\label{subsec:loss_function}

The KAN parameters $\boldsymbol{\theta}$ are determined by minimizing a composite physics-informed loss consisting of the interior PDE residual and the two clamped boundary residuals. A discrete least-squares formulation is given by
\begin{equation}
    \mathcal{L}(\boldsymbol{\theta})
    =
    \lambda_{\mathrm{PDE}}\mathcal{L}_{\mathrm{PDE}}
    +
    \lambda_u\mathcal{L}_{u}
    +
    \lambda_n\mathcal{L}_{n},
    \label{eq:total_loss}
\end{equation}
where $\lambda_{\mathrm{PDE}}$, $\lambda_u$, and $\lambda_n$ are nonnegative weighting parameters.

The interior physics loss is
\begin{equation}
    \mathcal{L}_{\mathrm{PDE}}
    =
    \frac{1}{N_{\Omega}}
    \sum_{k=1}^{N_{\Omega}}
    \left|
    r_{\mathrm{PDE}}
    (\mathbf{x}_k;\boldsymbol{\theta})
    \right|^2.
    \label{eq:pde_loss}
\end{equation}

The displacement boundary loss is
\begin{equation}
    \mathcal{L}_{u}
    =
    \frac{1}{N_{\partial\Omega}}
    \sum_{k=1}^{N_{\partial\Omega}}
    \left|
    r_u
    (\mathbf{x}_k^b;\boldsymbol{\theta})
    \right|^2,
    \label{eq:u_loss}
\end{equation}
and the normal-slope boundary loss is
\begin{equation}
    \mathcal{L}_{n}
    =
    \frac{1}{N_{\partial\Omega}}
    \sum_{k=1}^{N_{\partial\Omega}}
    \left|
    r_n
    (\mathbf{x}_k^b;\boldsymbol{\theta})
    \right|^2.
    \label{eq:n_loss}
\end{equation}

Consequently, the trained KAN-PINN approximation is defined by
\begin{equation}
    \boldsymbol{\theta}^{*}
    =
    \underset{\boldsymbol{\theta}}{\operatorname{arg\,min}}
    \;
    \mathcal{L}(\boldsymbol{\theta}),
    \label{eq:optimization}
\end{equation}
and the resulting approximation to the plate displacement is
\begin{equation}
    u_{\mathrm{KAN}}(x,y)
    =
    u_{\boldsymbol{\theta}^{*}}(x,y).
    \label{eq:kan_solution}
\end{equation}

\subsection{Continuous least-squares interpretation}
\label{subsec:continuous_loss}

The discrete loss in \eqref{eq:total_loss} can be interpreted as a quadrature approximation to a continuous residual functional. In particular, one may define
\begin{equation}
\begin{aligned}
    \mathcal{J}(u)
    &=
    \lambda_{\mathrm{PDE}}
    \left\|
    \Delta^2u-f
    \right\|_{L^2(\Omega)}^2
    \\
    &\quad+
    \lambda_u
    \left\|
    u-g
    \right\|_{L^2(\partial\Omega)}^2
    +
    \lambda_n
    \left\|
    \frac{\partial u}{\partial n}-j
    \right\|_{L^2(\partial\Omega)}^2.
\end{aligned}
\label{eq:continuous_loss}
\end{equation}

The exact solution of the boundary-value problem satisfies $\mathcal{J}(u)=0$, provided that all residuals vanish. Thus, the KAN-PINN seeks a function within the chosen KAN hypothesis space for which the interior differential equation and the clamped boundary conditions are simultaneously satisfied.

\section{Neural network architecture}
\label{sec:architecture}

Unlike conventional mesh-based methods or standard Multi-Layer Perceptrons (MLPs) with fixed activation functions, our framework represents the unknown transverse deflection $u(x,y)$ globally and mesh-free via a Kolmogorov--Arnold Network (KAN) embedded within a physics-informed training pipeline. As illustrated in Figure~\ref{fig:kan_architecture}, the KAN architecture replaces traditional node-based activations with learnable 1-D functions placed directly on the network's edges.

\subsection{Kolmogorov--Arnold network layer formulation}
In an ordinary neural network, each connection between nodes is parameterized by a single scalar weight. In contrast, every connection in a KAN is a small, learnable 1-D function consisting of a linear transformation combined with a sum of basis functions (such as radial basis functions or B-splines). Specifically, for a layer mapping an input vector $\mathbf{x} \in \mathbb{R}^{i_{\text{in}}}$ to an output vector $\mathbf{o} \in \mathbb{R}^{o_{\text{out}}}$, each output component $out_o$ is computed as:
\begin{equation}
    out_o = b_o + \sum_{i} \left[ w_{oi} x_i + \sum_{k=1}^{n_b} c_{oik} \exp\left( - \left( \frac{x_i - \mu_{ik}}{\sigma_{ik}} \right)^2 \right) \right],
\end{equation}
where $b_o$ is a learnable bias, $w_{oi}$ represents the linear weight (residual connection), and $c_{oik}$ denotes the learnable coefficient of the $k$-th Gaussian bump. 

\begin{figure}[H]
    \centering
    \includegraphics[width=0.95\textwidth]{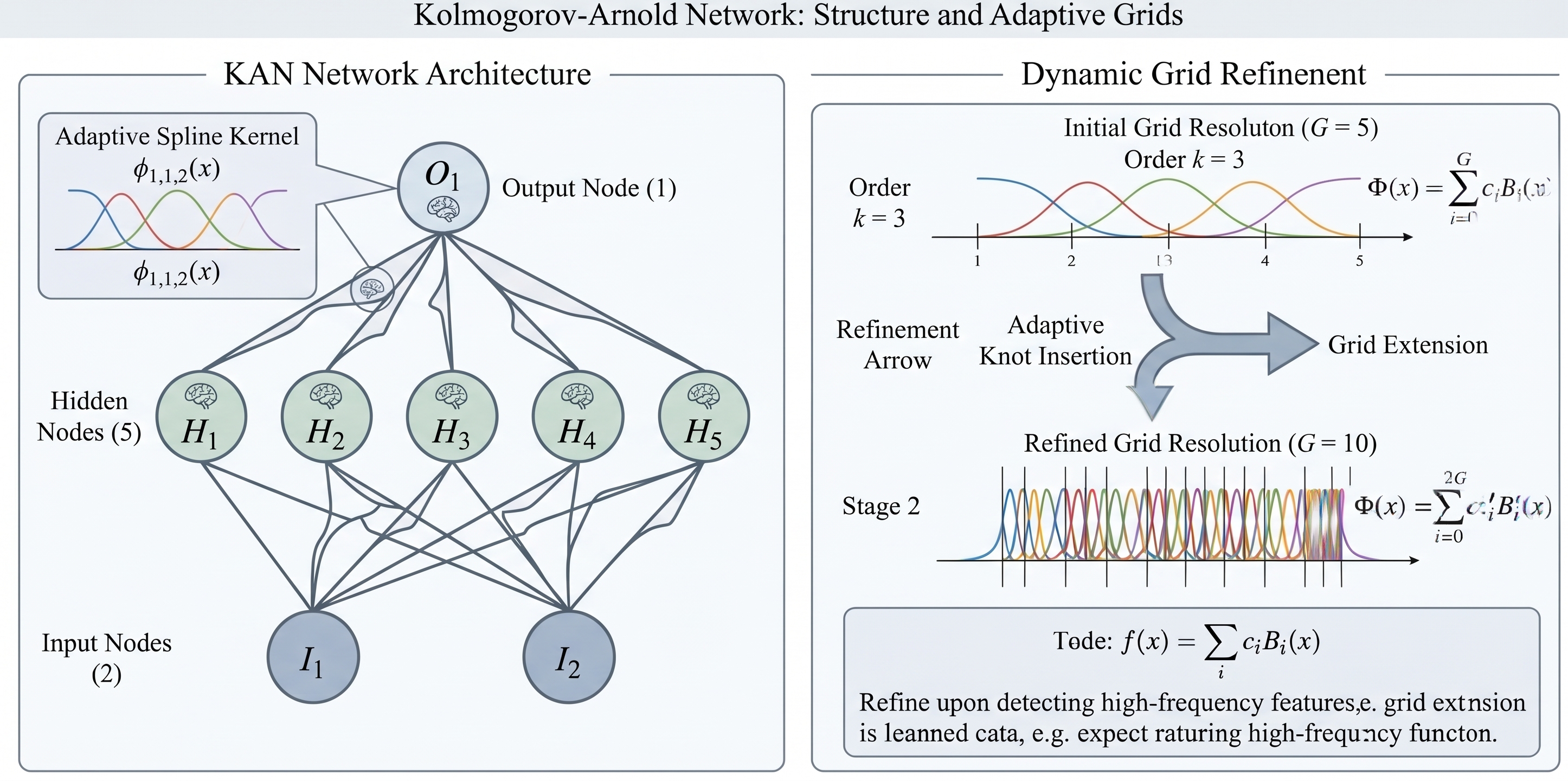}
    \caption{Kolmogorov--Arnold Network structure featuring learnable 1-D functions along the edges between input, hidden, and output nodes (left), and the dynamic grid refinement process utilizing adaptive knot insertion and grid extension (right).}
    \label{fig:kan_architecture}
\end{figure}

To ensure that high-order derivatives remain well-defined and reproducible, the Gaussian bump centers $\mu_{ik}$ are arranged on a fixed grid of $n_b$ points spread uniformly across each layer's input range, with widths $\sigma_{ik}$ fixed proportionally to the center spacing. Furthermore, as depicted in the right panel of Figure~\ref{fig:kan_architecture}, the framework utilizes \textit{dynamic grid refinement}. This mechanism allows the network to adaptively refine its resolution from an initial coarse grid to a finer grid via adaptive knot insertion and grid extension upon detecting localized high-frequency features or sharp gradients.

While Multi-Layer Perceptrons (MLPs) assign fixed nonlinear activation functions to nodes and linear weights to edges, Kolmogorov-Arnold Networks (KANs) invert this paradigm by placing learnable univariate functions on the edges and performing simple linear summation at the nodes. Given a supervised learning problem with input-output pairs $\{x_i, y_i\}$, KANs operationalize the Kolmogorov-Arnold representation theorem by parameterizing each individual edge function as a B-spline curve.

\subsection{Approximation guarantees and neural scaling behavior}
The expressive power of deep KANs is governed by rigorous approximation theory that extends the classical Kolmogorov-Arnold representation to smooth, arbitrarily wide and deep architectures. Under the assumption that a target function $f(x)$ admits a representation through $(k+1)$-times continuously differentiable layer functions $\Phi_l$, the approximation error bounded by the grid size $G$ satisfies a strict convergence rate:
$$\|f - (\Phi_{L-1}^G \circ \dots \circ \Phi_0^G) x \|_{C^m} \le C G^{-k-1+m}$$
where $C$ is a constant independent of the grid resolution $G$. For standard cubic splines ($k=3$), this yields a theoretical scaling exponent of $\alpha = 4$.

\subsection{Activation function justification}
A critical requirement for fourth-order PINNs is the choice of activation function. Because the governing biharmonic operator $\Delta^2 u$ differentiates the network output four times consecutively, the activation function must possess smooth, non-vanishing high-order derivatives. Standard piecewise-linear activations such as ReLU have a zero second derivative almost everywhere, which would reduce the biharmonic residual to zero or nonsense. Consequently, the hyperbolic tangent ($\tanh$) activation function is chosen deliberately between layers to ensure continuous high-order differentiability.

\subsection{Network topology and parameterization}
The optimal network architecture is structured as a four-layer composition with a uniform hidden width of $20$, mapping spatial coordinates $(x,y)$ to the scalar deflection $u_\theta$:
\begin{equation}
    h_1 = \tanh(K_1(x,y)), \quad h_2 = \tanh(K_2 h_1), \quad h_3 = \tanh(K_3 h_2), \quad u_\theta = K_4 h_3,
\end{equation}
where $K_1, K_2, K_3,$ and $K_4$ denote the respective KAN layers. Table~\ref{tab:kan_layers} breaks down the layer dimensions, input ranges, and parameter counts, totaling a compact $7,801$ trainable parameters.

\begin{table}[htbp]
    \centering
    \caption{Layer-wise breakdown of the $2 \to 20 \to 20 \to 20 \to 1$ KAN architecture.}
    \label{tab:kan_layers}
    \begin{tabular}{lccccc}
        \hline
        \textbf{Layer} & \textbf{Mapping} & \textbf{Centre Range} & \textbf{Bumps / Edge} & \textbf{Activation} & \textbf{Parameters} \\ \hline
        $K_1$ & $\mathbb{R}^2 \to \mathbb{R}^{20}$ & $[0, 1]$ & 8 & $\tanh$ & 380 \\
        $K_2$ & $\mathbb{R}^{20} \to \mathbb{R}^{20}$ & $[-1, 1]$ & 8 & $\tanh$ & 3,620 \\
        $K_3$ & $\mathbb{R}^{20} \to \mathbb{R}^{20}$ & $[-1, 1]$ & 8 & $\tanh$ & 3,620 \\
        $K_4$ & $\mathbb{R}^{20} \to \mathbb{R}^1$ & $[-1, 1]$ & 8 & Linear & 181 \\ \hline
        \multicolumn{5}{l}{\textbf{Total Trainable Parameters}} & \textbf{7,801} \\ \hline
    \end{tabular}
\end{table}

\section{Optimization strategy and training procedure}
\label{sec:optimization}

Training physics-informed neural networks governed by fourth-order differential operators presents severe numerical stiffness. Each successive differentiation amplifies high-frequency noise inherent in parameter updates, causing standard gradient descent trajectories to stall or oscillate violently when minimizing raw fourth-order residuals. To overcome these optimization hurdles, we implement a robust multi-stage hybrid training pipeline coupled with an automated hyperparameter exploration loop.

\subsection{Hybrid first- and second-order optimization schedule}
Because gradient-based optimization of high-order PINNs requires both global basin exploration and rapid local convergence, training proceeds via a structured hybrid schedule:
\begin{enumerate}
    \item \textbf{Exploration phase (AdamW):} The optimization process initializes with the AdamW stochastic optimizer, acting as a global explorer to rapidly guide KAN parameters into the correct low-loss basin.
    \item \textbf{Precision polish phase (restarted L-BFGS):} Once AdamW plateaus, optimization transitions to a second-order quasi-Newton L-BFGS optimizer, executing an escalating, restarted polish schedule to descend rapidly through stiff ravines in the loss landscape.
\end{enumerate}

\subsection{Automated train-change-improve exploration loop}
Rather than relying on manual tuning, the optimization framework is driven by an automated $24/7$ hill-climbing search protocol. This automated protocol executes repeated training runs and applies structural or hyperparameter adjustments whenever optimization plateaus. 

\begin{algorithm}[H]
\caption{Automated train-change-improve search and hybrid optimization pipeline}
\label{alg:optimization_loop}
\begin{algorithmic}[1]
\State Initialize KAN parameters $\boldsymbol{\theta}$, collocation sets $\mathcal{X}_{\Omega}$ and $\mathcal{X}_{\partial\Omega}$.
\State Set search space for penalty weights ($\lambda_{\text{dir}}, \lambda_{\text{neu}}$) and optimization schedules.
\For{$iteration = 1$ to $124$}
    \State \textbf{Phase 1 (Exploration):} Minimize composite loss $\mathcal{L}(\boldsymbol{\theta})$ using AdamW for $N_{\text{adam}}$ epochs.
    \State \textbf{Phase 2 (Polish):} Apply restarted L-BFGS optimizer for precision convergence.
    \State Evaluate validation loss and compute ground-truth relative $L_2$ error against $u_{\text{exact}}$.
    \If{current model improves upon previous best ($19$ successful improvements total)}
        \State Save optimal checkpoint $\boldsymbol{\theta}^*$ and log hyperparameters.
    \Else
        \State Apply hill-climbing adjustment (e.g., scale $\lambda_{\text{dir}}$, modify L-BFGS restart frequency).
    \EndIf
\EndFor
\State \textbf{Return} Best optimized model variant .
\end{algorithmic}
\end{algorithm}

Across an extensive search comprising $124$ automated training attempts, the loop achieved $19$ successive performance improvements. The final optimal configuration (, in \texttt{float32} precision on CPU hardware) successfully balanced the boundary penalty weights ($\lambda_{\text{dir}} = 480$, $\lambda_{\text{neu}} = 30$).

\subsection{Manufactured solution and convergence monitoring}
To rigorously quantify accuracy, we employ the manufactured-solution technique on the unit square with an analytical target field:
\begin{equation}
    u_{\text{exact}}(x,y) = \sin(\pi x)\sin(\pi y).
\end{equation}
Applying the biharmonic operator analytically yields the exact source term:
\begin{equation}
    f(x,y) = \Delta^2 u_{\text{exact}}(x,y) = 4\pi^4 \sin(\pi x)\sin(\pi y).
\end{equation}
Training monitors both the composite loss $\mathcal{L}$ and the ground-truth relative $L_2$ error against $u_{\text{exact}}$, driving training loss down to $3.065 \times 10^{-6}$ and relative $L_2$ error to $1.593 \times 10^{-5}$ ($0.00159\%$).

\section{Numerical results and error analysis}
\label{sec:results}

To evaluate the predictive accuracy and computational efficacy of the proposed mesh-free KAN-PINN framework, we analyze the performance of the best-performing model variant trained using the direct fourth-order formulation on CPU hardware with \texttt{float32} precision.

\subsection{Quantitative accuracy metrics}
The quantitative performance of the final optimized model is summarized in Table~\ref{tab:numerical_results}. Despite utilizing a compact architecture consisting of only $7,801$ trainable parameters, the network achieves exceptional precision in approximating the manufactured biharmonic solution $u_{\text{exact}}(x,y) = \sin(\pi x)\sin(\pi y)$.

\begin{table}[H]
    \centering
    \caption{Quantitative performance metrics for the optimal biharmonic KAN-PINN model .}
    \label{tab:numerical_results}
    \begin{tabular}{lc}
        \hline
        \textbf{Metric} & \textbf{Value} \\ \hline
        Final Training Loss ($\mathcal{L}$) & $3.065 \times 10^{-6}$ \\
        Relative $L_2$ Error & $1.593 \times 10^{-5}$ ($0.00159\%$) \\
        Mean Absolute Error & $5.935 \times 10^{-6}$ \\
        Maximum Pointwise Error & $3.256 \times 10^{-5}$ \\
        Trainable Parameters & $7,801$ \\
        Precision \& Hardware & \texttt{float32} CPU \\ \hline
    \end{tabular}
\end{table}

\subsection{Field approximations, pointwise errors, and residual distributions}
\begin{figure}[H]
    \centering
    \includegraphics[width=0.98\textwidth]{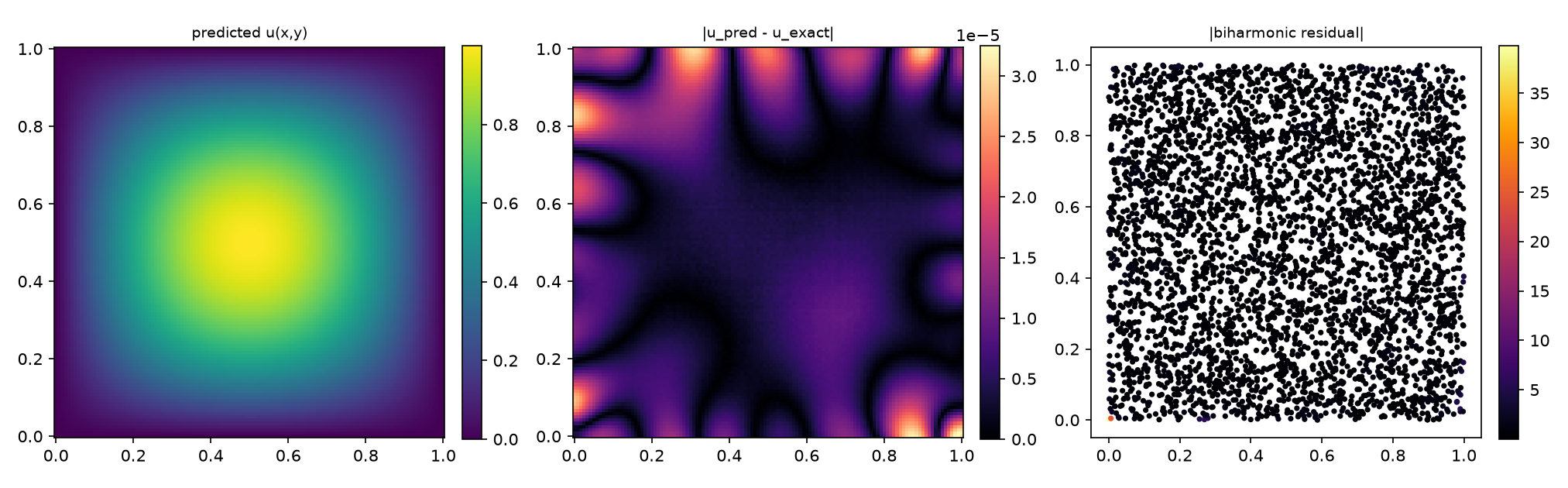}
    \caption{Spatial distribution of the predicted deflection $u(x,y)$ (left), absolute pointwise error $|u_{\text{pred}} - u_{\text{exact}}|$ (middle), and absolute biharmonic residual $|r_{\text{PDE}}|$ across collocation points (right) for model.}
    \label{fig:solution_fields}
\end{figure}

Figure~\ref{fig:solution_fields} illustrates the spatial distributions of the predicted transverse deflection $u(x,y)$, the absolute pointwise error magnitude $|u_{\text{pred}} - u_{\text{exact}}|$, and the absolute biharmonic residual field $|r_{\text{PDE}}|$ across the computational domain $\Omega$. 

As shown in the left panel, the predicted deflection profile accurately captures the smooth sinusoidal behavior characteristic of the clamped plate benchmark without artificial smoothing or mesh artifacts. The absolute error distribution (middle panel) indicates that maximum discrepancies remain strictly bounded on the order of $10^{-5}$, with minor boundary layer variations attributable to the stringent fourth-order clamped boundary constraints. The right panel displays the absolute biharmonic residual evaluated across the interior collocation points, demonstrating uniform adherence to the governing differential equation throughout the domain.

\begin{figure}[H]
    \centering
    \includegraphics[width=0.5\textwidth]{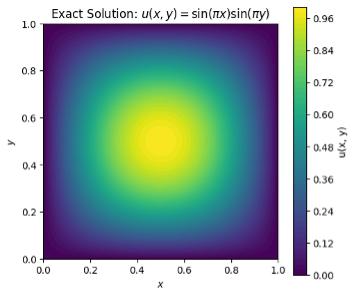}
    \caption{Analytical manufactured exact solution $u_{\text{exact}}(x,y) = \sin(\pi x)\sin(\pi y)$ on the unit square domain.}
    \label{fig:exact_solution}
\end{figure}

To visually validate the accuracy of the optimal model, Figure~\ref{fig:exact_solution} presents the analytical manufactured solution $u_{\text{exact}}(x,y) = \sin(\pi x)\sin(\pi y)$ across the unit square domain $\Omega$. 

When comparing the analytical solution in Figure~\ref{fig:exact_solution} directly with the network's predicted deflection field $u_{\text{pred}}(x,y)$ (displayed in the left panel of Figure~\ref{fig:solution_fields}), the KAN-PINN architecture demonstrates exceptional qualitative agreement. Both fields exhibit a smooth, symmetric peak value of $1.0$ at the spatial center $(0.5, 0.5)$ and rigorously satisfy the homogeneous clamped boundary conditions ($u = 0$) along all four edges of the domain. Furthermore, the absence of boundary layer leakage, artificial numerical diffusion, or Gibbs oscillations confirms that the edge-based B-spline parameterization accurately captures both the global amplitude and local curvature of the fourth-order biharmonic operator without requiring a physical grid mesh.

\subsection{Training and validation loss dynamics}
\begin{figure}[H]
    \centering
    \includegraphics[width=0.85\textwidth]{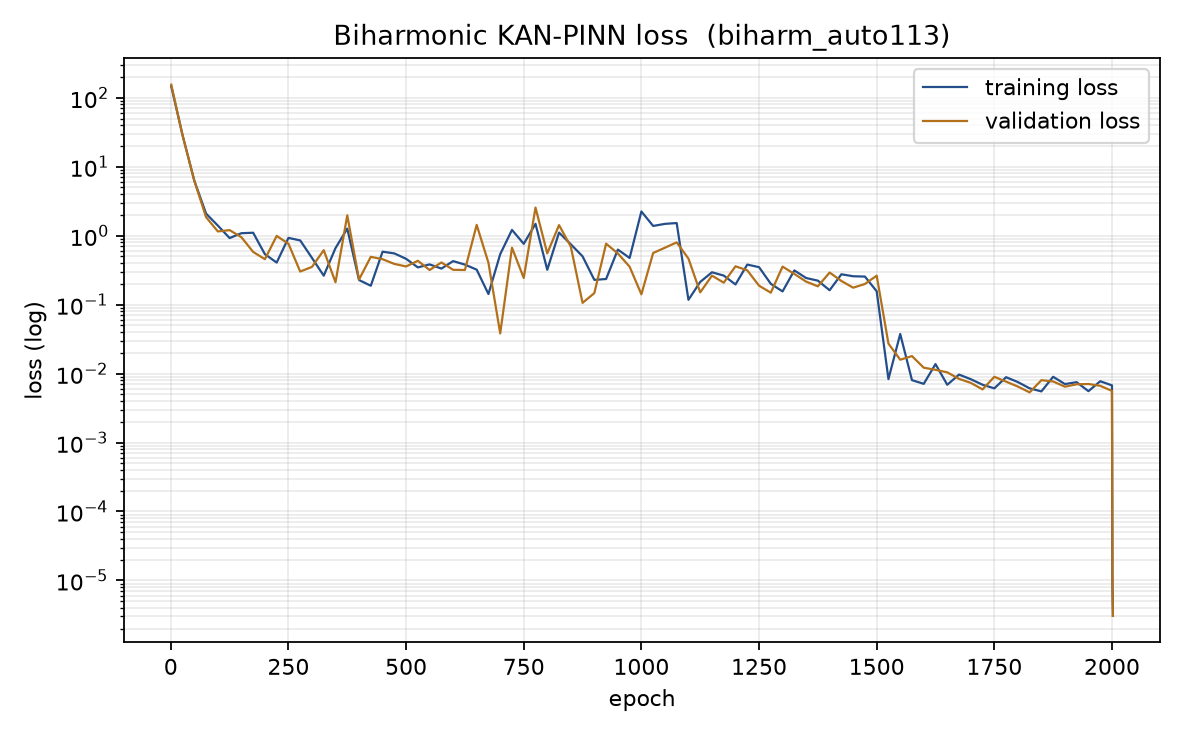}
    \caption{Training and validation loss history over $2,000$ epochs for the optimal biharmonic KAN-PINN model , highlighting the sharp convergence phase driven by the restarted L-BFGS optimizer.}
    \label{fig:loss_history}
\end{figure}

The temporal evolution of the training and validation losses over $2,000$ epochs for the model is presented in Figure~\ref{fig:loss_history}. The loss curve reflects the multi-stage hybrid optimization strategy: an initial exploration phase characterized by steady descent under AdamW, followed by a sharp order-of-magnitude drop in loss as the restarted L-BFGS polish schedule takes effect in the final optimization epochs. Both training and validation loss curves track each other closely, confirming robust generalization across the collocation points without overfitting.

\subsection{Error reduction via automated exploration}
Throughout the automated $24/7$ hill-climbing search comprising $124$ total training attempts, the algorithm successfully captured $19$ successive performance improvements. The automated train-change-improve loop systematically drove the error down from an initial baseline of $0.24\%$ to a final relative $L_2$ error of $0.00159\%$ ($1.593 \times 10^{-5}$), representing an overall error reduction factor exceeding $150\times$. 

\begin{figure}[H]
    \centering
    \includegraphics[width=0.88\textwidth]{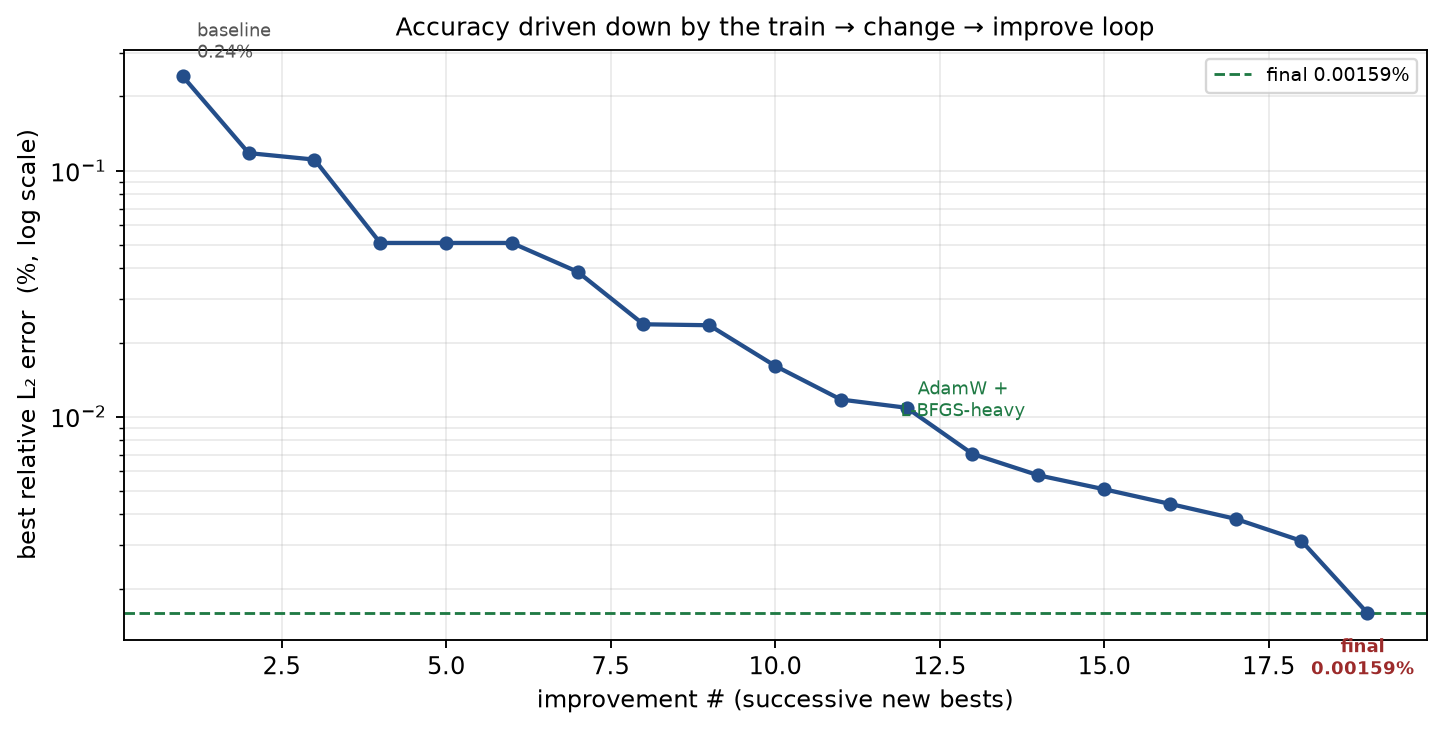}
    \caption{Progression of the best relative $L_2$ error (log scale) across the $19$ successive automated improvement steps during the train-change-improve search loop, culminating in the final error of $0.00159\%$.}
    \label{fig:accuracy_trajectory}
\end{figure}

Figure~\ref{fig:accuracy_trajectory} illustrates this multi-step convergence trajectory across the successive new best models. As highlighted by the search trajectory, the decisive levers governing this high degree of convergence were the escalating, restarted L-BFGS polish schedule and the implementation of a strong Dirichlet boundary penalty weight ($\lambda_{\text{dir}} = 480$), which effectively controlled the severe numerical stiffness associated with the fourth-order biharmonic operator boundaries. Conversely, tests involving larger network widths, double-precision (\texttt{float64}), and increased collocation sampling densities yielded negligible performance gains, confirming that optimization dynamics—rather than over-parameterization—dominate fourth-order PINN training.

\section{Conclusion}
\label{sec:conclusion}

A mesh-free numerical framework based on Kolmogorov--Arnold Physics-Informed Neural Networks (KAN-PINNs) has been successfully established and evaluated for the high-accuracy approximation of fourth-order boundary value problems governed by the biharmonic equation. By integrating learnable univariate functions parameterized via radial basis functions on network edges, smooth hyperbolic tangent activations for robust high-order differentiability, a direct normalized residual formulation, and a hybrid AdamW-to-L-BFGS optimization schedule, the severe numerical stiffness associated with fully clamped plate boundary conditions was effectively surmounted.

Through an automated hyperparameter exploration loop, critical training levers—namely, the restarted L-BFGS polish schedule and rigorous Dirichlet penalty calibration—were systematically optimized. The resultant compact architecture, comprising only $7,801$ trainable parameters, attained a final relative $L_2$ error of $1.593 \times 10^{-5}$ ($0.00159\%$) and a training loss of $3.065 \times 10^{-6}$ on the manufactured benchmark. These findings demonstrate that Kolmogorov--Arnold networks provide a powerful, mesh-free alternative to traditional multilayer perceptrons and classical discretization techniques for resolving complex, high-order differential equations without requiring explicit mesh generation or auxiliary variable transformations.

Future extensions of this adaptive KAN-PINN paradigm will focus on investigating complex irregular geometries, three-dimensional plate and shell structures, coupled multi-physics systems, and transient dynamic vibration phenomena in scientific and engineering computations.

\bibliographystyle{plain}  
\bibliography{KAN_Biharmonic}

\end{document}